\documentclass[reqno]{amsart}
\usepackage{enumerate}
\usepackage{amsthm}
\usepackage[colorlinks=true]{hyperref}
\usepackage[all]{xypic}
\usepackage{bbm}
\usepackage[final]{epsfig}
\usepackage{psfrag}

\theoremstyle{theorem}
\newtheorem{theorem}{Theorem}
\newtheorem{proposition}{Proposition}
\newtheorem{corollary}{Corollary}
\newtheorem{lemma}{Lemma}

\theoremstyle{definition}
\newtheorem{definition}{Definition}

\theoremstyle{remark}
\newtheorem{remark}{Remark}

\newtheorem{remarksAux}{Remarks}

\newcommand{\defin}[1]{\emph{#1}}
\newcommand{\R}{\mathbb{R}}
\newcommand{\N}{\mathbb{N}}
\newcommand{\Z}{\mathbb{Z}}

\newcommand{\epsi}{\varepsilon}
\newcommand{\om}{\omega}

\newcommand{\fr}{\texttt{\textit{fr}}}
\newcommand{\oc}{\texttt{\textit{oc}}}

\title{Coding Discretizations of Continuous Functions}

\author{Cristobal Rojas}
\author{Serge Troubetzkoy}
\begin{document}

\maketitle

\begin{abstract}
We consider several coding discretizations of continuous functions which reflect their variation at some given precision. We study certain statistical and combinatorial properties of the sequence of finite words obtained by coding a typical continuous function when the diameter of the discretization tends to zero.  Our main result is that any finite word appears on a subsequence discretization with any desired limit frequency. 
\end{abstract}

\section{Introduction}

Take a straight line in the plane and code it by a $0-1$ sequence as
follows: each time it crosses an integer vertical line (that is, $x=n$ for some $n\in \Z$) write a $0$ and each
time it crosses an integer horizontal line ($y=n$ for some $n\in \Z$) write a $1$.  In the case of
irrational slope the corresponding sequence is
called a Sturmian sequence \cite{F}. A classical result tells us that each word
that
appears in such a sequence has a limiting frequency. Moreover, the set of numbers occurring as limit frequencies can be completely described  \cite{VB}.
Recently similar condings have been considered for quadratic functions
and limiting frequencies are caluculated for words which appear
\cite{DTZ}. 

In this article we ask the question if limiting frequencies can appear
in more general circumstance: namely for typical, in the sense of
Baire, continuous functions. For such functions it is not clear which
kind of coding should be used.  Here we propose three different notions
of coding.  For each of these codings we study two different questions:
if all finite words can appear in a code or not, and
if words in the code of a typical function can have a limiting
frequency.

A \defin{discretization system} of $[0,1]$ is a sequence $X_{n}:=\{0=x^n_{1},x^n_{2},...,x^n_{N_{n}}=1\} \subset [0,1]$ where, 
\begin{enumerate}
\item $X_{1}\subset X_{2} \subset ... \subset X_{n} \subset .... \subset [0,1]$,
\item For each $X_{n}$, $x^n_{i}<x^n_{i+1}$ for all $1\leq i < N_{n}$,
\item The \defin{maximal resolution} $H_{n}:=\max_{1 \leq i < N_{n}} |x^n_{i+1}-x^n_{i}|$ converges to zero.
\end{enumerate}
We denote by $h_{n}:=\min_{1\leq i < N_{n}} |x^n_{i+1}-x^n_{i}|$, the \defin{minimal resolution}. To each discretization system $X_{n}$, we  associate the (uniform) discretization of the image space 
given by $Y_{n}:=\{y^n_{j}=jh_{n}:j\in \N\}$. 

\vspace{.2cm}
Let $f\in C([0,1])$.  For each $x_{i}^n \in X_{n}$ there is a unique $j\in \Z$
such that $f(x_{i}^n)\in[y_{j}^n,y_{(j+1)}^n)$. 
Let us denote this $j$ by $f_{i}^n$.

\begin{figure}\label{fig1}
\psfrag*{x2}{\footnotesize$x_2$}
\psfrag*{x3}{\footnotesize$x_3$}
\psfrag*{x4}{\footnotesize$x_4$}
\psfrag*{x5}{\footnotesize$x_5$}
\psfrag*{x6}{\footnotesize$x_6$}
\psfrag*{x7}{\footnotesize$x_7$}
\psfrag*{0}{\footnotesize$0=x_1$}
\psfrag*{1}{\footnotesize$1=x_8$}
\psfig{file=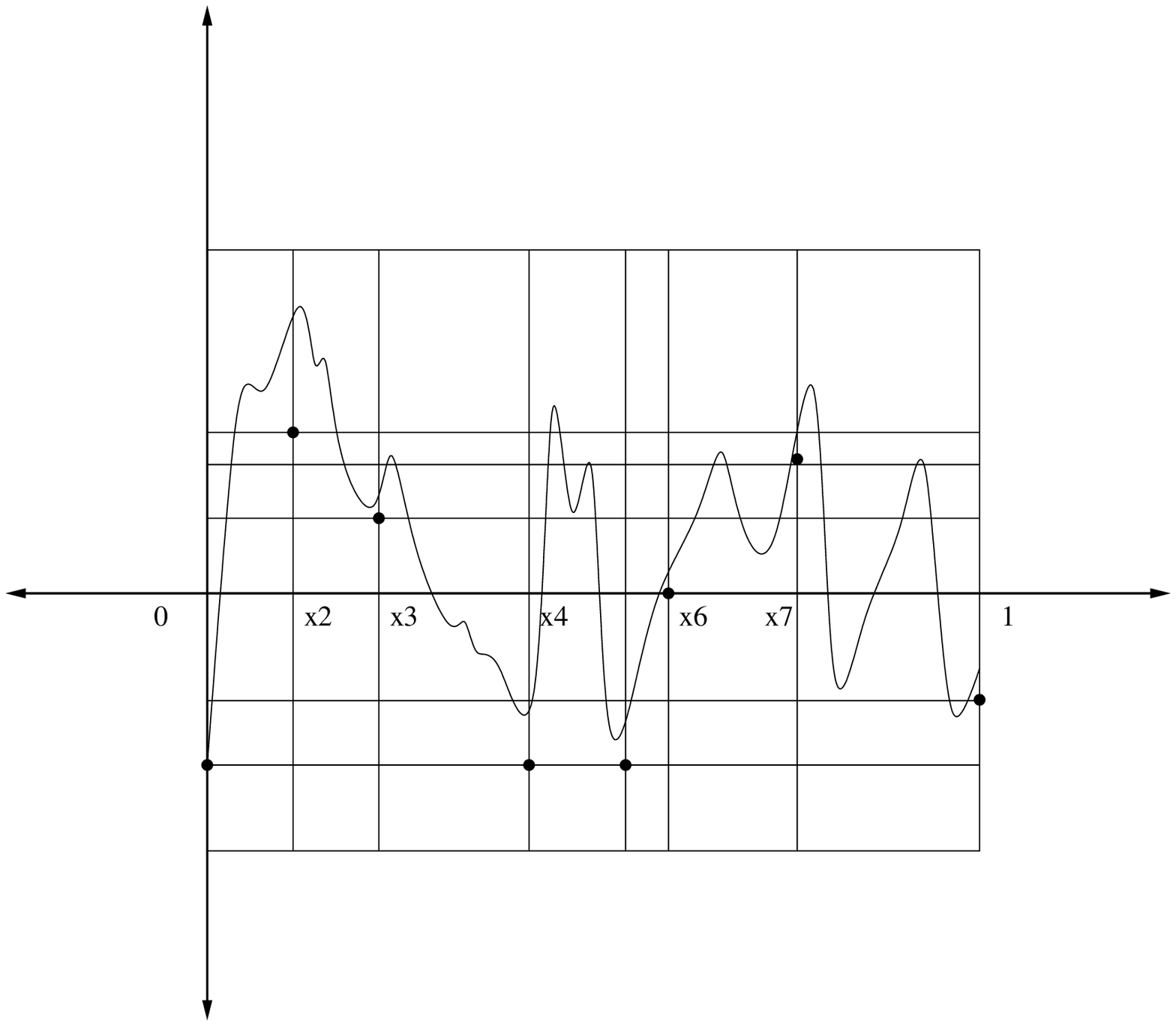,height=100mm}
\caption{The various codes considered: \newline
quantitative: $Q(f,n) =5^-2^-3 0 2 2^-3$ \newline
qualitative: $q(f,n) = 1^-1^-1 0 1 1^-1$ \newline
stretched: $s(f,n) = 111110^-1^-10^-1^-1^-100110110^-1^-1^-10$}
\end{figure}

\begin{definition} (See Figure \ref{fig1}) The \defin{quantitative code} $Q(f,n)\in \Z^{N_{n}-1}$ of $f\in C([0,1])$ is defined by: 

\[
	Q(f,n)_{i}:= f_{i+2}^n-f_{i+1}^n, \phantom{qqq} 0\leq i<N_{n}-1.	
\]

The  \defin{qualitative} version $q(f,n)\in \{-1,0,1\}^{N_{n}-1}$ of the quantitative code $Q(f,n)$ is defined by  setting

\[
	q(f,n)_{i}:=
		\begin{cases}
  			 \phantom{-}1,  &\mbox{ if }  Q(f,n)_{i}>0  \\
			 \phantom{-}0,  &\mbox{ if }  Q(f,n)_{i}=0 \\
  			 		  -1,  &\mbox{ if } Q(f,n)_{i}<0.
	\end{cases}
\]

Finally, the \defin{stretched} version $s(f,n)\in \{1,0,-1\}^*$ of the quantitative code
$Q(f,n)$ is defined as follows:  if $Q(f,n)_i$ is positive then we replace it by a run of $Q(f,n)_i$ 1's
followed by a zero and by a run of $-1$'s followed by a zero if $Q(f,n)_i$ is negative.
\end{definition}

All three of these codes  seem natural in terms of discrete curves on the
computer screen.  In case when the discretization system is uniform, the streched quantitative
code of a line segment with irrational slope is exactly the well known coding by Sturmian sequences
\cite{F}.

Let us introduce some more notation in order to state our main results.
Let $w,v$ be finite words over the same alphabet $\Sigma$ (finite or infintie) such that $|w|\leq |v|$. We denote by
\begin{equation}
\oc(w,v):=\#\{j: v_{j}^{j+|w|}=w, 0\leq j \leq |v|-|w|\}
\end{equation}
the number of times $w$ occurs in $v$ and by $\fr(w,v):=\frac{\oc(w,v)}{|v|}$  the relative frequency of $w$ in $v$. The \defin{minimal periodic factor length} $p(w)$ of $w$ is defined to be $p(w):=\min\{|u|:\oc(w,wu)=2\}$. For example, $p(010)=2$.

\begin{remark}The limit relative frequency  of $w$ in some infinite sequence $v^n$ is at most $\frac{1}{p(w)}$. That is:
$$
\limsup_{n}\fr(w,v^n)\leq\frac{1}{p(w)}.
$$
\end{remark}

Our main result is the following:

\begin{theorem}\label{t}Let $X_{n}$ be a discretization system. For a typical $f\in C([0,1])$ the following holds:
\begin{enumerate}[\upshape(i)]
\item\label{p:quali} (Qualitative) For any  $w \in \{-1,0,1\}^*$ and $\alpha\in [0,\frac{1}{p(w)}]$, there exists a subsequence $n_{i}$ such that 
$$
\lim_{i\to \infty}\fr(w,q(f,n_{i}))=\alpha.
$$

\item\label{p:quanti} (Quantitative) Suppose that $X_{n}$ satisfies $\liminf_{n}nh_{n}=0$. Then for any  $w \in \Z^*$ and $\alpha\in [0,\frac{1}{p(w)}]$, there exists a subsequence $n_{i}$ such that 
$$
\lim_{i\to \infty}\fr(w,Q(f,n_{i}))=\alpha.
$$

\item\label{p:st} (Stretched)  Suppose that $X_{n}$ satisfies $\liminf_{n}nh_{n}=0$ and $\frac{H_{n}}{h_{n}}$ is bounded. Then 
$$
\liminf_{n\to \infty}\fr(0,s(f,n))=0,
$$
if $f(1)\geq f(0)$

$$
\liminf_{n\to \infty}\fr(1,s(f,n))=\limsup_{n\to \infty}\fr(-1,s(f,n))=\frac{1}{2}, 
$$

\noindent and if $f(1)\leq f(0)$
$$
\limsup_{n\to \infty}\fr(1,s(f,n))=\liminf_{n\to \infty}\fr(-1,s(f,n))=\frac{1}{2} 
$$

\end{enumerate}

\end{theorem}



\section{Preliminaries}

We start by a simple result, which says that one can focus on functions which do not intersect the discretization.

\begin{lemma}\label{l.Cd}
Let $X_{n}$ be a discretization system. Then for a typical function $f$ one has that for all $n\in \N$ and all $i=1,...,N_{n}$, $f(x^n_{i})\in (y_j^{n},y^n_{j+1})$, for the corresponding $j\in \N$.
\end{lemma}
\begin{proof}
The set $F_{n}=\{f: f(x^n_{i})\neq y_j^{n} \text{ for all }j\in \N \text{ and }i=1,...,N_{n}\}$ is clearly open and dense. Hence, $\bigcap_{n}F_{n}$
is a $G_{\delta}$-dense set.
\end{proof}

One would expect that codings of  typical functions contains  few zeros and all possible words of $1$'s and $-1's$.  This is partially true.

\begin{proposition}\label{no.zeros}
Let $X_{n}$ be a discretization system.  For a typical $f$,   $q(f,n)$ contains no $0$, infinitely often.
\end{proposition}
\begin{proof}
 We prove that the set of functions such that for all $n\in \N$, there exists $m\geq n$ such that  $q(f,m)_{i}\neq0$ for all $i=0,...,N_{m}-2$, is residual in $C([0,1])$. Observe that $q(f,n)_{i}\neq 0$ whenever $|f(x^n_{i+1})-f(x^n_{i})|>h_{n}$.  Clearly, the set 
$$
F^m=\{f: |f(x^m_{i+1})-f(x^m_{i})|>h_{m} \text{ for all } i=1,...,N_{m}-1\}
$$
is open. Moreover, for each $n\in \N$, the set 
$$
\bigcup_{m\geq n}F^m
$$

is a dense open set. Indeed, given $g\in C([0,1])$ and $\epsi>0$,  there exists $m\geq n$ such that $h_{m}< \epsi$ and it is easy to construct a function  $f\in F^{m}$ such that $\Vert g-f\Vert_{\infty}<\epsi$. Therefore, 
$$
\bigcap_{n}\bigcup_{m\geq n}F^m
$$
is a $G_{\delta}$-dense set.

\end{proof}

\begin{remark} In the previous result, the symbol $0$ cannot be replaced by $1$ nor by $-1$.  On the other hand, Theorem \ref{t}  says that the qualitative and the quantitative codings of a typical function does not posses any statistical regularity. So that from a statistical viewpoint, the symbols $1$ or $-1$ (or any $n\in \Z$ in the quantitative case) are not privileged with respect to $0$. 
\end{remark}

\subsection{Approximation by  $\epsi$-boxes}

Here we will describe a simple construction which will be used in the proofs of our main results.

 Let $\delta>0$. For a given $n$ we define a subdiscretization $X_{n}^{\delta}:=\{x_{i_{k}}:k=1,...,K\}$ of $X_{n}$ as follows:
\[
\begin{array}{rl}
x_{i_{1}}=&0,\\
x_{i_{k+1}}=&{\max\{x_{i}^n\in X_{n}:x_{i}^n<x_{i_{k}}+2\delta\}}\\
x_{i_{K}}=&1
\end{array}
\]
The number of points of $X_{n}$ in the interval $(x_{i_{k}},x_{i_{k+1}}]$ will be denoted by $l_{k}$. With this notation we have $x_{i_{k+1}}=x_{i_{k}+l_{k}}$. 

\begin{figure}\label{fig2}
\centerline{\psfig{file=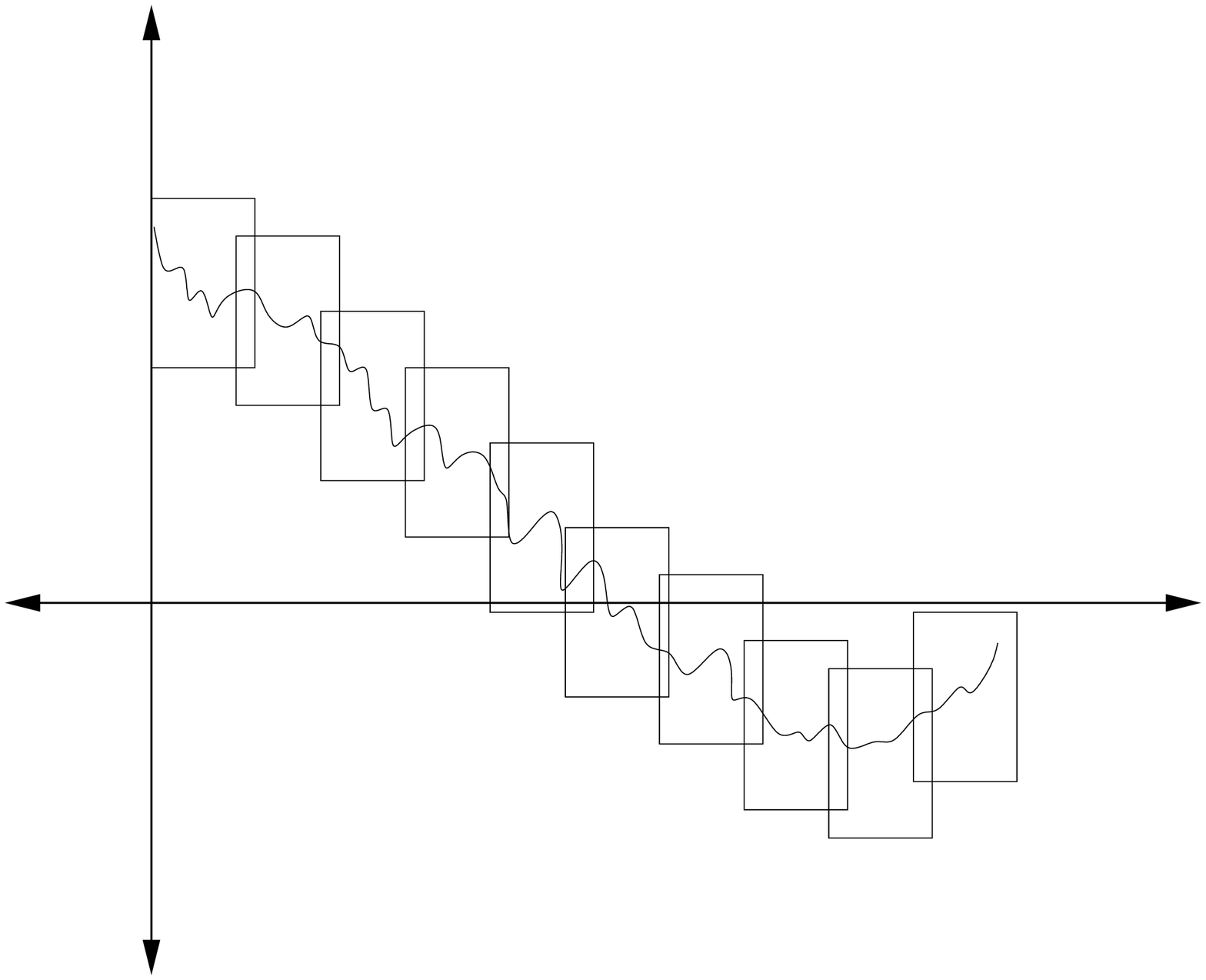,height=100mm}}
\caption{An $\epsi$-cover by the boxes $B_k(g,\epsi,\delta)$}
\end{figure}

Next, to each $g\in C([0,1])$ and $\epsi>0$, the associated \defin{$\epsi$-boxes} $B_{k}(g,\epsi,\delta)$ are defined by:

\begin{equation}
B_{k}(g,\epsi,\delta):=(x_{i_{k}},x_{i_{k}}+2\delta)\times(g(\Delta_{k})-\frac{\epsi}{2},g(\Delta_{k})+\frac{\epsi}{2}))
\end{equation}
where $\Delta_{k}=\frac{x_{i_{k}}+x_{i_{k+1}}}{2}$. See figure
\ref{fig2}. 
We shall write just $B_{k}$ when no confusion is possible.

Let $\delta_{g}:\R^+\to\R^+$ denote the modulus of continuity of $g$. That is, for every  $x,x'$ in $[0,1]$, if $|x-x'|<\delta_{g}(\epsi)$ then $|f(x)-f(x')|<\epsi$.  

\begin{lemma}\label{boxes} For the  $\epsi$-boxes $B_{k}(g,\epsi,\delta)$, $k=1,...,K$, the following holds:
\begin{enumerate}[\upshape(i)]
\item If $\delta<\delta_{g}(\frac{\epsi}{2})$ then the $\epsi$-boxes  form an $\epsi$-cover of the graph of $g$. That is,  any $(x,y)\in \bigcup_{k}B_{k}(g,\epsi,\delta_{g}(\frac{\epsi}{2}))$ satisfy $|g(x)-y|<\epsi$.  
\item If $\lceil\frac{1}{2\delta}\rceil H_{n}<2\delta$, then $K= \lceil\frac{1}{2\delta}\rceil+1$. 
\end{enumerate}
\end{lemma}
\begin{proof}
Let $(x,y)\in \bigcup_{k}B_{k}$.  Then $(x,y)\in B_{k}$ for some $k$. Hence $|x-\Delta_{k}|<\delta<\delta_{g}(\frac{\epsi}{2})$ which implies $|g(x)-g(\Delta_{k})|<\frac{\epsi}{2}$.  Since $g(\Delta_{k})-\frac{\epsi}{2}<y<g(\Delta_{k})+\frac{\epsi}{2}$ we conclude $|g(x)-y|<\epsi$.

 Since $x_{i_{k+1}}<x_{i_{k}}+2\delta$, we have that $K\geq \lceil\frac{1}{2\delta}\rceil+1$.  Now, for each $k$ we have $x_{i_{k}}+2\delta - x_{i_{k+1}}\leq H_{n}$. It follows that 
$$
K\leq  \left\lceil\frac{1}{2\delta}\right \rceil +  \left\lceil \frac{\lceil \frac{1}{2\delta}\rceil H_{n}}{2\delta}\right \rceil.
$$
Hence, if  $2\delta>\lceil\frac{1}{2\delta}\rceil H_{n}$ we obtain $K\leq \lceil\frac{1}{2\delta}\rceil+1$. 
\end{proof}

\subsection{Words and frequencies}Consider a finite word $w$ over some alphabet $\Sigma$. For each $\alpha \in [0,\frac{1}{p(w)}]$ and $t>0$, it is easy to construct a sequence of finite words $v_{k}$, $k=1,...,K-1$, satisfying $|v_{k}|=l_{k}-1$ and $|\fr(w,v_{k})-\alpha|<\frac{1}{3t}$.  Let $b_{k}\in \Sigma$ be any sequence of $K-1$ letters and put $v=v_{1}b_{1}v_{2}b_{2}\cdot \cdot \cdot v_{K-1}b_{K-1}$. 

\begin{lemma}\label{freq}
If $\frac{(K-1)|w|}{|v|}<\frac{1}{3t}$ then $|\fr(w,v)-\alpha|<\frac{1}{t}$.
\end{lemma}
\begin{proof}
 Put $\oc(w,v_{k})=p_{k}$, then we have
\begin{equation}\label{1}
\frac{\sum_{k=1}^{K-1} p_{k}}{|v|}\leq \fr(w,f,n)\leq \frac{\sum_{k=1}^{K-1} p_{k}}{|v|}+\frac{(K-1)|w|}{|v|}
\end{equation}

 A simple calculation yields
\begin{equation}\label{2}
\frac{\sum_{k=1}^{K-1} p_{k}}{|v|}=\frac{\sum_{k=1}^{K-1}p_{k}}{\sum_{k=1}^{K-1}(l_{k}-1)} - \frac{\sum_{k=1}^{K-1}p_{k}}{(\sum_{k=1}^{K-1}(l_{k}-1))^2+(K-1)\sum_{k=1}^{K-1}(l_{k}-1)}.
\end{equation}

On the one hand we have:

$$
\left\vert\frac{\sum_{k=1}^{K-1}p_{k}}{\sum_{k=1}^{K-1}(l_{k}-1)}-\alpha \right\vert=\left|\frac{1}{K-2}\sum_{k=1}^{K-1}\fr(w,v_{k}) - \alpha\right|<\frac{1}{3t}
$$
and on the other hand, the absolute value of the second term in the right side of equation (\ref{2}) is less than
$$
\frac{|v|-(K-1)}{(|v|-(K-1))^2+(|v|-(K-1))(K-1)}=\frac{1}{|v|}\leq\frac{1}{3t}
$$

so that 
$$
\left|\frac{\sum_{k=1}^{K-1} p_{k}}{|v|}-\alpha\right|<\frac{2}{3t}.
$$
Since  $\frac{(K-1)|w|}{|v|}<\frac{1}{3t}$, from equation (\ref{1})  we obtain $|\fr(w,f,n)-\alpha|\leq\frac{1}{t}$ and the lemma is proved.
\end{proof}


\section{Proofs}

\begin{proof}[Proof of Theorem \ref{t}]  We begin by proving parts
(\ref{p:quali}) and (\ref{p:quanti}). 
For each finite word $w$ in $\{-1,0,1\}^*$ or $\Z^*$,  let $\{\alpha^w_{s}\}_{s\in \N}$ be a sequence which is dense on $[0,\frac{1}{p(w)}]$. Let $F_{i}$ denote  the open sets  defined in Lemma \ref{l.Cd} (the set of functions which do not intersect the discretization $X_{i}$).  For integers $s,n,t$, consider the sets

$$
\overline{F}^q_{w,s,n,t}:=\{f\in \cap_{i\leq n}F_{i}: |\fr(w,q(f,n))-\alpha_{s}|\leq\frac{1}{t}\},
$$

$$
\overline{F}^Q_{w,s,n,t}:=\{f\in  \cap_{i\leq n}F_{i}:  |\fr(w,q(f,n))-\alpha_{s}|\leq\frac{1}{t}\}.
$$

Clearly these sets are open since a function $f$ in $\cap_{i\leq n}F_{i}$ can be perturbed without changing  its code $q(f,n)$ or $Q(f,n)$. Hence, the following sets are open too.

$$
F^q_{w,s,m,t}:=\{f: \exists n\geq m, f\in \overline{F}^q_{w,s,n,t} \},
$$

$$
F^Q_{w,s,m,t}:=\{f: \exists n\geq m,  f\in \overline{F}^Q_{w,s,n,t} \}.
$$

We now show that these sets are moreover dense. Let $g\in C([0,1])$ and $\epsi>0$.  We will construct a function $f$ in $F^q_{w,s,m,t}$ (respectively $F^Q_{w,s,m,t}$) such that $\Vert f-g \Vert_{\infty}\leq \epsi$. 

\vspace{.2cm}
\noindent \textbf{Case $\boldsymbol{F^q_{w,s,m,t}}$}. Put $\delta<\min\{\delta_{g}(\frac{\epsi}{2}),\frac{\epsi}{4}\}$ and  let $B_{k}$ be the associated $\epsi$-boxes.
For $n\geq m$ large enough (in particular such that $\lceil\frac{1}{2\delta}\rceil H_{n}<2\delta$) there exists a sequence of finite words  $v_{k}$, $k=1,...,K-1$,  such that $|v_{k}|=l_{k}-1$, $|\fr(w,v_{k})-\alpha_{s}|<\frac{1}{3t}$ and $\frac{(K-1)|w|}{N_{n}}<\frac{1}{3t}$.

\begin{figure}\label{fig3}
\psfrag*{xik}{\footnotesize{$x_{i_k}$}}
\psfrag*{xik1}{\footnotesize{$x_{i_{k}+1}$}}
\psfrag*{xik2}{\footnotesize{$x_{i_{k+1}}$}}
\psfrag*{a}{\footnotesize{$a_k$}}
\psfrag*{b}{\footnotesize{$b_k$}}
\psfrag*{c}{\footnotesize{$c_k$}}
\centerline{\psfig{file=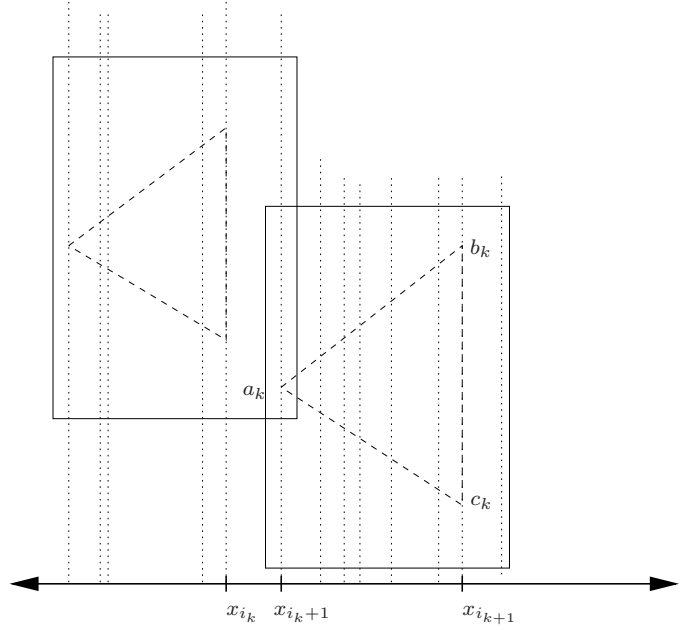,height=100mm}}
\caption{The triangles $(a_k,b_k,c_k).$}
\end{figure}

We claim that a function $f_{v}$ can be constructed such that for each $k$ we have $q(f_{v},n)_{i_{k}+1}^{i_{k}+l_{k}}=v_{k}$ and $f_{v}$ is $\epsi$-close to $g$ (the interval $(x_{i_{k}},x_{i_{k}+1})$ is reserved to make ``the bridge'' and there are $K$ such intervals, see figure 3). To see this, observe that the condition $2\delta<\frac{\epsi}{2}$ implies that for each $k$, the triangles of vertices $(a_{k},b_{k},c_{k})$ defined by 
\begin{eqnarray*}
a^q_{k} &= &(x_{i_{k}+1},g(\Delta_{k}))\\
b^q_{k} &= &(x_{i_{k+1}}, g(\Delta_{k})+|v_{k}|h_{n} )\\
c^q_{k} &= &(x_{i_{k+1}}, g(\Delta_{k})-|v_{k}|h_{n} )
\end{eqnarray*}
 are included in $B_{k}$ and that for any $v_{k}$, a function $f_{v}$ such that $q(f_{v},n)_{i_{k}+1}^{i_{k}+l_{k}}=v_{k}$ can be inscribed in these triangles. By lemma \ref{boxes} a function $f$ so constructed satisfies $\Vert f-g \Vert_{\infty}\leq \epsi$. By lemma \ref{freq} we have that $|\fr(w,q(f,n))-\alpha_{s}|<\frac{1}{t}$.

\vspace{.2cm}
\noindent \textbf{Case $\boldsymbol{F^Q_{w,s,m,t}}$}. The proof that these sets are dense is the same as for the sets $F^q_{w,s,m,t}$, with the only exception that we have to take $\delta<\frac{\epsi}{4H(w)}$ where $H(w):=\max_{i}|w_{i}|$ denotes the \defin{hight} of $w$. This condition assures that  a function $f_{v}$ such that $q(f_{v},n)_{i_{k}+1}^{i_{k}+l_{k}}=v_{k}$ can be inscribed in the corresponding triangles:

 \begin{eqnarray*}
a^Q_{k} &= &(x_{i_{k}+1},g(\Delta_{k}))\\
b^Q_{k} &= &(x_{i_{k+1}}, g(\Delta_{k})+ |v_{k}|H(w)h_{n} )\\
c^Q_{k} &= &(x_{i_{k+1}}, g(\Delta_{k})-|v_{k}|H(w)h_{n} ).
\end{eqnarray*}

It follows that the sets
$$
\bigcap_{w,s,m,t}F^q_{w,s,m,t}     \phantom{qqq} \text{ and } \phantom{qq} \bigcap_{w,s,m,t}F^Q_{w,s,m,t} 
$$
are both $G_{\delta}$-dense.

Finally we prove part (\ref{p:st}).
Let  $u_{n}:=\oc(1,s(f,n))$ be the numbers of  $1$'s (or ``ups'') in $s(f,n)$ and $d_{n}:=\oc(-1,s(f,n))$ be the number of  $-1$'s (or ``downs''),  in stage $n$. Then $V_{n}=u_{n}+d_{n}$ denotes the  \defin{total $n$-variation}. By definition we have that $|s(f,n)|=V_{n}+N_{n}$, where $N_{n}$ is both the cardinality of the discretization and the number of zeros. Hence we have
$$
\fr(1,s(f,n))= \frac{u_{n}}{V_{n}+N_{n}}.
$$

We will need the following lemma:

\begin{lemma}\label{l.quantitative}
Let $X_{n}$ be a discretization system satisfying $\liminf_{n}nh_{n}= 0$. Then, for a typical $f$, there are infinitely many $n$ such that $ u_{n}>\frac{nN_{n}}{3}$ and $d_{n}>\frac{nN_{n}}{3}$.
\end{lemma}
\begin{proof}
Consider the set of functions 
$$
\overline{F}_{n}:=\left\{f: card\{i: Q(f,n)_{i}>n\}>\frac{N_{n}}{3} \text{ and }  card\{i: Q(f,n)_{i}<-n\}>\frac{N_{n}}{3}\right\}\cap \bigcap_{i\leq n}F_{i}
$$
 This is an open set. Moreover, for any $m\in \N$, the set
$$
\bigcup_{n\geq m}\overline{F}_{n}
$$
is dense. For let $g\in C[0,1]$ and consider the associated $\epsi$-boxes $B_{k}$. It is clear that for some  $n\geq m$ such that $nh_{n}<\frac{\epsi}{2}$ one can construct a function $f$ satisfying $graph(f)\subset \cup_{k}B_{k}$ and $|f(x_{i+1})-f(x_{i})|>n$ for all $i$.  Moreover, we can alternate the sign of  $|f(x_{i+1})-f(x_{i})|$ at every $i$, with at most $K$ exceptions. Hence the function so constructed belongs to $\overline{F}_{n}$ and then the set 
$$
\bigcap_{m}\bigcup_{n\geq m}\overline{F}_{n}
$$
is $G_{\delta}$ dense. 

\end{proof}

Now, a simple calculation yields
$$
 \frac{u_{n}}{V_{n}}=\frac{1}{2-\frac{\Delta_{n}}{u_{n}}}
$$
where $\Delta_{n}=\lfloor  \frac{(f(1)-f(0))}{h_{n}} \rfloor$. So, if $f(1)=f(0)$ we have $\frac{u_{n}}{V_{n}}=\frac{1}{2}$.
 Let $M$ be a bound for  $\frac{H_{n}}{h_{n}}$. We have then that $\frac{1}{h_{n}}\leq MN_{n}$ and hence $\Delta_{n}\leq (f(1)-f(0))MN_{n}$.  By  Lemma \ref{l.quantitative} we have that

$$
\frac{(f(1)-f(0))MN_{n}}{u_{n}}<\frac{3(f(1)-f(0))MN_{n}}{nN_{n}}
$$
and 
$$
\frac{N_{n}}{V_{n}}<\frac{3}{n}
$$
for infinitely many $n$, so that, 
\begin{eqnarray}
\liminf \frac{\Delta_{n}}{u_{n}}=0 & \text{ if } & f(1)>f(0), \\
\limsup  \frac{\Delta_{n}}{u_{n}}=0 & \text{ if }& f(1)<f(0).
\end{eqnarray}

Hence, when $f(1)>f(0)$ we have
$$
\liminf_{n\to \infty}\frac{u_{n}}{V_{n}+N_{n}} =   \liminf_{n\to \infty} \frac{u_{n}}{V_{n}}=\frac{1}{2-\liminf \frac{\Delta_{n}}{u_{n}}}=\frac{1}{2}
$$
and  when $f(1)<f(0)$ we have
$$
\limsup_{n\to \infty}\frac{u_{n}}{V_{n}+N_{n}} = \limsup_{n\to\infty} \frac{u_{n}}{V_{n}}=\frac{1}{2-\liminf \frac{\Delta_{n}}{u_{n}}}=\frac{1}{2}
$$
and the results follows by symmetry.

\end{proof}




\begin{thebibliography}{99}

\bibitem[B]{VB}   	
V. Berth\'e, {\em Frequencies of Sturmian series factors},
Theoretical Computer Science. 165(2): 295-309 (1996)

\bibitem[DTZ]{DTZ}A.\ Daurat, M.\ Tajine, M.\ Zouaoui,  {\em Fr\'equences des motifs d'une discr\'etisation
    de courbe}\\
 http://www.lama.univ-savoie.fr/gdrim-geodis/images/21novembre/daurat.pdf

\bibitem[F]{F}   	
{\em Substitutions in dynamics, arithmetics, and combinatorics},
N.\ Pythease Fogg, Springer Lecture Notes 1794.
Chapter 6. Sturmian Sequences, by P.\ Arnoux


\end{thebibliography}
\end{document}